\documentclass[12pt]{article}
\usepackage[letterpaper, left=1.5in, right=1in,top=1in,bottom=1in]{geometry}
\usepackage{authblk}
\usepackage{tikz}
\usepackage{chngcntr}
\counterwithin{figure}{section}
\usetikzlibrary{arrows}
\usepackage{placeins}	
\usepackage{amsfonts}
\usepackage{graphicx}
\usepackage{standalone}
\usepackage{amsmath, amsthm, amssymb}
\usepackage[capitalize, noabbrev]{cleveref}
\newtheorem{thrm}{Theorem}[section]
\newtheorem{prop}[thrm]{Proposition}

\newtheorem{conj}[thrm]{Conjecture}

\newenvironment{pf}           {\noindent{\bf Proof:} }%
                                {\null\hfill$\Box$\par\medskip}
                           
\usepackage[nottoc]{tocbibind}

\begin{document}

\title{Independent Coverings and Orthogonal Colourings}
\author[ ]{Kyle MacKeigan}
\affil[ ]{Department of Mathematics, Dalhousie University}
\affil[ ]{Email: \textit{Kyle.m.mackeigan@gmail.com}}
\maketitle

\begin{abstract}
In this paper, two open conjectures are disproved. One conjecture regards independent coverings of sparse partite graphs, whereas the other conjecture regards orthogonal colourings of tree graphs. A relation between independent coverings and orthogonal colourings is established. This relation is applied to find independent coverings of some sparse partite graphs. Additionally, a degree condition providing the existence of an independent covering in the case where the graph has a square number of vertices is found.
\end{abstract}

\section{Introduction}
Two vertex colourings of a graph are \textit{orthogonal} if they have the property that when two vertices are coloured with the same colour in one colouring, then those vertices receive different colours in the other colouring. An \textit{orthogonal colouring} of a graph is a pair of orthogonal vertex colourings. The \textit{orthogonal chromatic number} of a graph $G$, denoted by $O\chi(G)$, is the minimum number of colours required for a proper orthogonal colouring.

It will be shown that orthogonal colourings generalize independent coverings, which are now defined. An \textit{independent transversal} of a graph, with respect to a vertex partition $P$, is an independent set that contains exactly one vertex from each vertex class. An \textit{independent covering} of a graph, with respect to a vertex partition $P$, is a collection of disjoint independent transversals with respect to $P$ that spans all of the vertices. 

For example, consider the rook graph, which is the graph that represents all legal moves of the rook chess piece on a chessboard. Note that a pair of orthogonal Latin squares of size $n$ gives an independent covering of this graph. This is because one Latin square provides a partition of the graph into independent sets, and the second orthogonal Latin square provides the independent transversals with respect to this partition.

Of special interest is the case where all parts of the vertex partition have the same size. In part, this is due to their applications to fields such as equitable colourings \cite{furmanczyk2006equitable}. This gives rise to the following definition. A graph is \textit{$[n,k,r]$-partite} if the vertices can be partitioned into $n$ independent sets of size $k$, with exactly $r$ independent edges between every pair of independent sets.

Note that in an $[n,k,r]$-partite graph, the parameter $n$ is the number of parts in a partition $P$ of the vertices into independent sets. Thus, this corresponds to a vertex colouring using $n$ colours. Then, since each part has size $k$, if an independent cover with respect to $P$ exists, then this corresponds to a second orthogonal colouring, using $k$ colours. In general, an independent covering with respect to a partition of the vertices into independent sets gives an orthogonal colouring. 

On the other hand, an orthogonal colouring only gives an independent covering if the sizes of the colour classes in the first colouring are the same and the sizes of the colour classes in the second colouring are the number of colours used in the first colouring. Therefore, independent coverings can be thought of as a special case of orthogonal colourings. Previously studied research is now discussed.

Orthogonal colourings were first defined in 1985 by Archdeacon, Dinitz, and Harary in the context of edge colourings \cite{archdeacon1985orthogonal}. Later in 1999, Caro and Yuster studied orthogonal colourings, this time in the context of vertex colourings \cite{caro1999orthogonal}. Then in 2013, Ballif studied upper bounds on sets of orthogonal vertex colourings \cite{ballif2013upper}. More recently, orthogonal colourings of Cayley graphs were explored \cite{janssen2020orthogonal}.

Some combinatorial objects, like list colourings, can be created by finding an independent transversal of a graph \cite{loh2007independent}. Thus, finding sufficient conditions for the existence of an independent transversal is important. This problem was originally studied by Bollob{\'a}s, Erd{\H{o}}s, and Szemer{\'e}di \cite{bollobas1975complete} in 1975. They conjectured that if each vertex class in a partition has size at least 2$\Delta$, where $\Delta$ is the maximum degree of the graph, then an independent transversal exists. This conjecture was later proved by Haxell \cite{haxell2001note} in 2001.

There is no degree condition for the existence of an independent covering. Thus, research on independent coverings is on specific families of graphs, like $[n,k,r]$-partite graphs. Independent transversals of $[n,k,r]$-partite graphs were originally studied by Erd{\H{o}}s, Gy{\'a}rf{\'a}s, and {\L}uczak \cite{erdHos1994independent} due to their usefulness in constructing hypergraphs with large girth and large chromatic number \cite{nevsetvril1979short}. For similar reasons, independent coverings of $[n,k,r]$-partite graphs were studied by Yuster \cite{yuster1997independent}.

Let $c(k,r)$ denote the maximal $n$ such that all $[n,k,r]$-partite graphs have an independent covering with respect to the given $[n,k,r]$-partition. Yuster \cite{yuster1997independent} showed that $k\geq c(k,r)\geq \min\{k,k-r+2\}$. Note that when $r=1,2$, the upper and lower bound coincide, giving that $c(k,1)=c(k,2)=k$. This led to the following conjecture:

\begin{conj}[\cite{yuster1997independent}]\label{Conjecture: Coverings}
For all $r\leq k$, $c(k,r)=k$.
\end{conj}

In this paper, Conjecture \ref{Conjecture: Coverings} is disproved by showing that there is a specific $[3,3,3]$-partite graph that does not have an independent covering with respect to the given $[3,3,3]$-partition. Then, a lower bound of $c(k,r)\geq \lceil\frac{k}{2}\rceil$ is obtained by using orthogonal colourings. This gives an improved lower bound on $c(k,r)$ for $r>\frac{k}{2}+2$. Lastly, a degree condition for independent coverings is obtained.

For a graph $G$ with $n$ vertices, $O\chi(G)\geq \lceil\sqrt{n}\,\rceil$ is an obvious lower bound. If $O\chi(G)=n$ and $G$ has $n^2$ vertices, then $G$ is said to have a \textit{perfect orthogonal colouring}. Note that in a perfect orthogonal colouring, all colour classes have the same size. Therefore, a degree condition for a graph to have a perfect orthogonal colourings gives a degree condition for a graph to have an independent covering. Such a condition is found by studying the following conjecture.

\begin{conj}[\cite{caro1999orthogonal}]\label{Conjecture: Trees}
If $T$ is a tree graph with $n$ vertices and $\Delta(T)<\frac{n}{2}$, then $O\chi(T)=\lceil\sqrt{n}\,\rceil$.
\end{conj}

Caro and Yuster proposed Conjecture \ref{Conjecture: Trees} because of a false categorization of the orthogonal chromatic number of double stars. The correct orthogonal chromatic number is determined in this paper and Conjecture \ref{Conjecture: Trees} is disproved. To conclude, we show that if $G$ is a $d$-degenerate graph with $n$ vertices and $\Delta(G)<\frac{\sqrt{n-2d-1}}{2}$, then $O\chi(G)=\lceil\sqrt{n}\,\rceil$.

\section{Independent Coverings of $[n,k,k]$-partite graphs}

Recall that $c(k,r)$ is the maximal $n$ such that all $[n,k,r]$-partite graphs have an independent covering with respect to the given $[n,k,r]$-partition. Conjecture \ref{Conjecture: Coverings} is disproved by showing that $c(3,3)\neq 3$. It will be shown later that an independent covering with respect to a different partition exists however.

\begin{thrm}\label{Theorem: Counter Example}
There exists a $[3,3,3]$-partite graphs that does not have an independent covering with respect to the $[3,3,3]$-partition. That is, $c(3,3)\neq 3$.
\end{thrm}
\begin{pf}
We will show that the $[3,3,3]$-partite graph $G$ in Figure \ref{Figure: Counter Example Graph} does not have an independent covering with respect to $P=\{\{x_0,x_1,x_2\},\{y_0,y_1,y_2\},\{z_0,z_1,z_2\}\}$. The three transversals are defined as $T_0$, $T_1$, $T_2$ and $T_0$, $T_1$, $T_2$ are populated later. Suppose for the sake of contradiction that $G$ does have an independent covering with respect to $P$ with the independent transversals $T_0, T_1,$ and $T_2$.  Without loss of generality, suppose that $x_0 \in T_0$, $x_1\in T_1$, and $x_2\in T_2$.

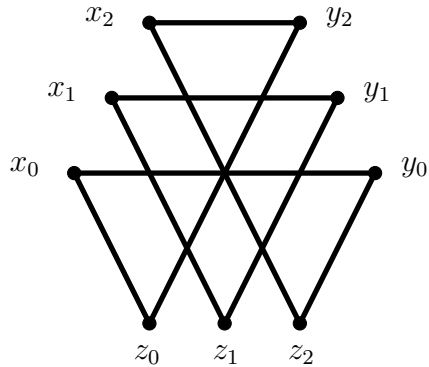
\begin{figure}[h!]
\centering
\begin{tikzpicture}[scale=1, line cap=round,line join=round,>=triangle 45,x=1.0cm,y=1.0cm]
\draw [line width=2.pt] (1.,4.)-- (3.,4.);
\draw [line width=2.pt] (0.5,3)-- (3.5,3);
\draw [line width=2.pt] (2.,0.)-- (3.5,3);
\draw [line width=2.pt] (0.5,3)-- (2.,0.);
\draw [line width=2.pt] (0.,2)-- (4.,2.);
\draw [line width=2.pt] (4.,2.)-- (3.,0.);
\draw [line width=2.pt] (3.,4.)-- (1.,0.);
\draw [line width=2.pt] (3,0)-- (1,4);
\draw [line width=2.pt] (1,0)-- (0,2);
\draw (0,4.35) node[anchor=north west] {$x_2$};
\draw (3.2,4.35) node[anchor=north west] {$y_2$};
\draw (0.65,-0.15) node[anchor=north west] {$z_0$};
\draw (-0.5,3.35) node[anchor=north west] {$x_1$};
\draw (3.7,3.35) node[anchor=north west] {$y_1$};
\draw (1.7,-0.15) node[anchor=north west] {$z_1$};
\draw (-1,2.35) node[anchor=north west] {$x_0$};
\draw (4.2,2.35) node[anchor=north west] {$y_0$};
\draw (2.7,-0.15) node[anchor=north west] {$z_2$};
\begin{scriptsize}
\draw [fill=black] (0.,2) circle (2.5pt);
\draw [fill=black] (0.5,3) circle (2.5pt);
\draw [fill=black] (1.,4.) circle (2.5pt);
\draw [fill=black] (3.,4.) circle (2.5pt);
\draw [fill=black] (3.5,3) circle (2.5pt);
\draw [fill=black] (4.,2.) circle (2.5pt);
\draw [fill=black] (1.,0.) circle (2.5pt);
\draw [fill=black] (2.,0.) circle (2.5pt);
\draw [fill=black] (3.,0.) circle (2.5pt);
\end{scriptsize}
\end{tikzpicture}
\caption{Counterexample Graph}
\label{Figure: Counter Example Graph}
\end{figure}

There are two properties used to show how the independent covering must be formed when $T_0$ is given. The first property is that each vertex in $\{x_0,x_1,x_2\}$, $\{y_0,y_1,y_2\}$, and $\{z_0,z_1,z_2\}$ must be in a different independent transversal. This is because an independent transversal can only have one vertex from each vertex class. The second property is that if two vertices are adjacent, then they must be in different transversals. This is because an independent transversal is an independent set.

\textbf{Case 1: $T_0=\{x_0,y_1,z_2\}$.} Since $x_2y_2\in E(G)$ and $x_2\in T_2$, $y_2\not\in T_2$ by the second property. Then, since $y_1\in T_0$, $y_2\not\in T_0$ by the first property. Therefore, $y_2\in T_1$ is the only option available. Now, if $z_1\in T_1$, then $x_1z_1\in E(G)$ and $x_1,z_1\in T_1$, which contradicts that $T_1$ is independent. If $z_0\in T_1$, then $y_2z_0\in E(G)$ and $y_2,z_0\in T_1$, which contradicts that $T_1$ is independent. 

\textbf{Case 2: $T_0=\{x_0,y_2,z_1\}$.} Since $x_1y_1\in E(G)$ and $x_1\in T_1$, $y_1\not\in T_1$ by the second property. Then, since $y_2\in T_0$, $y_1\not\in T_0$ by the first property. Therefore, $y_1\in T_2$ is the only option available. Thus, $y_0\in T_1$ is the only option available. Now, if $z_2\in T_2$, then $x_2z_2\in E(G)$ and $x_2,z_2\in T_2$, which contradicts that $T_2$ is independent. If $z_2\in T_1$, then $y_0z_2\in E(G)$ and $y_0,z_2\in T_1$, which contradicts that $T_1$ is independent.

\textbf{Case 3: $T_0=\{x_0,y_2,z_2\}$.} Since $x_1y_1\in E(G)$ and $x_1\in T_1$, $y_1\not\in T_1$ by the second property. Then, since $y_2\in T_0$, $y_1\not\in T_0$ by the first property. Therefore, $y_1\in T_2$ is the only option available. Now, if $z_1 \in T_{1}$, then $x_1z_1\in E(G)$ and $x_1,z_1\in T_1$, which contradicts that $T_1$ is independent. If $z_1\in T_2$, then $y_1z_1\in E(G)$ and $y_1,z_1\in T_2$, which contradicts that $T_2$ is independent.
\end{pf}

Theorem \ref{Theorem: Counter Example} illustrates that not every $[3,3,3]$-partite graphs has an independent covering with respect to the $[3,3,3]$-partition. The following theorem shows that all $[3,3,3]$-partite graphs have an independent covering with respect to some partition however. This is done by constructing an orthogonal colouring for each graph.

\begin{figure}[h!]
\centering
\begin{tikzpicture}[scale=1, line cap=round,line join=round,>=triangle 45,x=1.0cm,y=1.0cm]
\draw [line width=2.pt] (1.,4.)-- (3.,4.);
\draw [line width=2.pt] (0.5,3)-- (3.5,3);
\draw [line width=2.pt] (2.,0.)-- (3.5,3);
\draw [line width=2.pt] (0.5,3)-- (2.,0.);
\draw [line width=2.pt] (0.,2)-- (4.,2.);
\draw [line width=2.pt] (4.,2.)-- (3.,0.);
\draw [line width=2.pt] (3.,4.)-- (1.,0.);
\draw [line width=2.pt] (3,0)-- (1,4);
\draw [line width=2.pt] (1,0)-- (0,2);
\draw (-.35,4.35) node[anchor=north west] {$(0,0)$};
\draw (3.1,4.35) node[anchor=north west] {$(1,2)$};
\draw (0.25,-0.15) node[anchor=north west] {$(2,1)$};
\draw (-0.85,3.35) node[anchor=north west] {$(0,1)$};
\draw (3.6,3.35) node[anchor=north west] {$(1,0)$};
\draw (1.45,-0.15) node[anchor=north west] {$(2,2)$};
\draw (-1.35,2.35) node[anchor=north west] {$(0,2)$};
\draw (4.1,2.35) node[anchor=north west] {$(2,0)$};
\draw (2.65,-0.15) node[anchor=north west] {$(1,1)$};
\begin{scriptsize}
\draw [fill=black] (0.,2) circle (2.5pt);
\draw [fill=black] (0.5,3) circle (2.5pt);
\draw [fill=black] (1.,4.) circle (2.5pt);
\draw [fill=black] (3.,4.) circle (2.5pt);
\draw [fill=black] (3.5,3) circle (2.5pt);
\draw [fill=black] (4.,2.) circle (2.5pt);
\draw [fill=black] (1.,0.) circle (2.5pt);
\draw [fill=black] (2.,0.) circle (2.5pt);
\draw [fill=black] (3.,0.) circle (2.5pt);
\end{scriptsize}
\end{tikzpicture}
\caption{Orthogonal Colouring of $G_3$}
\label{Figure: G3 Colouring}
\end{figure}
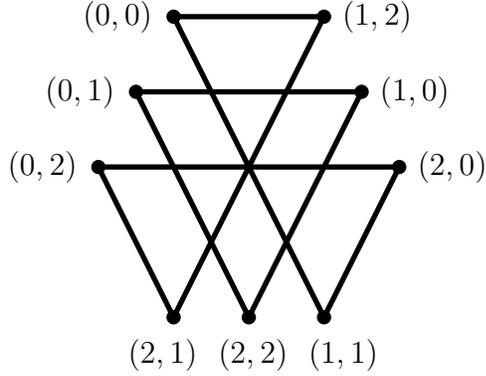

\begin{thrm}\label{Theorem: Independent Colouring}
Let $G$ be a $[3,3,3]$-partite graph, then $O\chi(G)=3$.
\end{thrm}
\begin{pf}
There are three different $[3,3,3]$-partite graphs. $G_1=K_3\cup K_3\cup K_3$, $G_2=C_9$, and $G_3=K_3\cup C_6$. Orthogonal colourings, where each colour class has the same size, of $G_1,G_2,$ and $G_3$ are given in Figure \ref{Figure: G1 Colouring}, Figure \ref{Figure: G2 Colouring}, and Figure \ref{Figure: G3 Colouring} respectively. The colour pairs represent the colours assigned by each of the two colourings.
\end{pf}

\begin{figure}[h!]
\centering
\begin{minipage}{.5\textwidth}
\centering
\begin{tikzpicture}[scale=.9, line cap=round,line join=round,>=triangle 45,x=1.0cm,y=1.0cm]
\draw [line width=2.pt] (1.,4.)-- (3.,4.);
\draw [line width=2.pt] (0.5,3)-- (3.5,3);
\draw [line width=2.pt] (2.,0.)-- (3.5,3);
\draw [line width=2.pt] (0.5,3)-- (2.,0.);
\draw [line width=2.pt] (0.,2)-- (4.,2.);
\draw [line width=2.pt] (4.,2.)-- (3.,0.);
\draw [line width=2.pt] (3.,4.)-- (1.,0.);
\draw [line width=2.pt] (1.,0.)-- (1.,4.);
\draw [line width=2.pt] (3.,0.)-- (0.,2);
\draw (-.35,4.35) node[anchor=north west] {$(0,0)$};
\draw (3.1,4.35) node[anchor=north west] {$(1,1)$};
\draw (0.25,-0.15) node[anchor=north west] {$(2,2)$};
\draw (-0.85,3.35) node[anchor=north west] {$(0,1)$};
\draw (3.6,3.35) node[anchor=north west] {$(1,2)$};
\draw (1.45,-0.15) node[anchor=north west] {$(2,0)$};
\draw (-1.35,2.35) node[anchor=north west] {$(0,2)$};
\draw (4.1,2.35) node[anchor=north west] {$(1,0)$};
\draw (2.65,-0.15) node[anchor=north west] {$(2,1)$};
\begin{scriptsize}
\draw [fill=black] (0.,2) circle (2.5pt);
\draw [fill=black] (0.5,3) circle (2.5pt);
\draw [fill=black] (1.,4.) circle (2.5pt);
\draw [fill=black] (3.,4.) circle (2.5pt);
\draw [fill=black] (3.5,3) circle (2.5pt);
\draw [fill=black] (4.,2.) circle (2.5pt);
\draw [fill=black] (1.,0.) circle (2.5pt);
\draw [fill=black] (2.,0.) circle (2.5pt);
\draw [fill=black] (3.,0.) circle (2.5pt);
\end{scriptsize}
\end{tikzpicture}
\caption{Orthogonal Colouring of $G_1$}
\label{Figure: G1 Colouring}
\end{minipage}%
\begin{minipage}{.5\textwidth}
\centering
\centering
\begin{tikzpicture}[scale=.9, line cap=round,line join=round,>=triangle 45,x=1.0cm,y=1.0cm]
\draw [line width=2.pt] (1.,4.)-- (3.,4.);
\draw [line width=2.pt] (0.5,3)-- (3.5,3);
\draw [line width=2.pt] (2.,0.)-- (3.5,3);
\draw [line width=2.pt] (0.5,3)-- (3,0);
\draw [line width=2.pt] (0.,2)-- (4.,2.);
\draw [line width=2.pt] (4.,2.)-- (3.,0.);
\draw [line width=2.pt] (3.,4.)-- (1.,0.);
\draw [line width=2.pt] (2,0)-- (1,4);
\draw [line width=2.pt] (1,0)-- (0,2);
\draw (-.35,4.35) node[anchor=north west] {$(0,0)$};
\draw (3.1,4.35) node[anchor=north west] {$(1,1)$};
\draw (0.25,-0.15) node[anchor=north west] {$(2,0)$};
\draw (-0.85,3.35) node[anchor=north west] {$(0,1)$};
\draw (3.6,3.35) node[anchor=north west] {$(1,2)$};
\draw (1.45,-0.15) node[anchor=north west] {$(2,1)$};
\draw (-1.35,2.35) node[anchor=north west] {$(0,2)$};
\draw (4.1,2.35) node[anchor=north west] {$(1,0)$};
\draw (2.65,-0.15) node[anchor=north west] {$(2,2)$};
\begin{scriptsize}
\draw [fill=black] (0.,2) circle (2.5pt);
\draw [fill=black] (0.5,3) circle (2.5pt);
\draw [fill=black] (1.,4.) circle (2.5pt);
\draw [fill=black] (3.,4.) circle (2.5pt);
\draw [fill=black] (3.5,3) circle (2.5pt);
\draw [fill=black] (4.,2.) circle (2.5pt);
\draw [fill=black] (1.,0.) circle (2.5pt);
\draw [fill=black] (2.,0.) circle (2.5pt);
\draw [fill=black] (3.,0.) circle (2.5pt);
\end{scriptsize}
\end{tikzpicture}
\caption{Orthogonal Colouring of $G_2$}
\label{Figure: G2 Colouring}
\end{minipage}
\end{figure}

Theorem \ref{Theorem: Counter Example} suggests that if $n>\lceil\frac{k}{2}\rceil$, then an independent covering with respect to the $[n,k,k]$-partition may not exist, but an independent covering with respect to some partition might. However, if $n\leq \lceil\frac{k}{2}\rceil$, then an independent covering with respect to the $[n,k,k]$-partition does exist. This is shown in the following theorem.

\begin{thrm}\label{Theorem: First Half}
All $[\frac{k}{2},k,k]$-partite graphs have an independent covering with respect to the $[\frac{k}{2},k,k]$-partition. That is, $c(k,k)\geq \lceil\frac{k}{2}\rceil$.
\end{thrm}
\begin{pf}
Let $n=\lceil\frac{k}{2}\rceil$ and let $G$ be an $[n,k,k]$-partite graph. To prove the statement, we need to show there is an independent covering with respect to the $[n,k,k]$-partition. We do this by constructing an orthogonal colouring such that the sizes of the colour classes in the first colouring are the same and the sizes of the colour classes in the second colouring are the number of colours used in the first colouring.  This way, the orthogonal colouring will correspond to an independent covering.

For $1\leq i \leq n$ and $1\leq j \leq k$, let $A_i$ denote the $i$-th vertex class and let $v_{i,j}$ denote the $j$-th vertex in $A_i$. For the first colouring, define $f_1(v_{i,j})=i$. Since each $A_i$ is an independent set, there is no colour conflict. For $1\leq m\leq n$, an independent covering with respect to the $[n,k,k]$-partition is established by constructing a second colouring $f_2$ using induction on $m$.  It will be shown that for each $i$, that each vertex of $A_i$ can be assigned a unique colour while not causing any colour conflicts.

For the base case $m=1$, define the second colouring on $A_1$ as $f_2(v_{1,j})=j$. Since each $v_{1,j}$ receives a different colour, there is no colour conflict. Suppose now that for $1< m \leq n$, that the second colouring is properly defined on the sets $A_1,\dots, A_{m-1}$ where each colour appears exactly once on each $A_i$. It is now shown that $f_2$ can be defined on $A_{m}$.

Let $C_j$ denote the set of available colours for the vertex $v_{k,j}$ in the second colouring. By definition, each vertex in $A_m$ has exactly one neighbour in each of $A_1,\dots, A_{m-1}$. Therefore, at most $m-1$ colours are appearing on the neighbours of each $v_{k,j}$. Since $f_2$ uses $k$ colours, this means there are at least $k-(m-1)\geq n$ colours available for $v_{k,j}$. Therefore, it follows that $|C_{j}|\geq n$. By the induction hypothesis, all colours appear exactly once on each $A_i$. Therefore, each colour appears in at least $k-(m-1)\geq n$ of the $C_{j}$. 

A corollary of Hall's Theorem says that if there are $n$ elements and $n$ sets where each set contains $k$ elements and each element appears in $k$ sets, then a system of distinct representatives exists \cite{van2001course}. Assigning the distinct representative of $C_j$ to $v_{k,j}$ will complete the orthogonal colouring. Since each colour class has the same size, this orthogonal colouring corresponds to an independent covering.
\end{pf}

It remains an open problem to show that this is the best possible lower bound. As seen in this section, there is a relation between independent coverings and orthogonal colourings of graphs with exactly $n^2$ vertices and $n$ colours. This relation will be explored further in the next section.

\section{Orthogonal Colourings of Tree Graphs}

In this section, a degree condition that guarantees the existence of an orthogonal colouring using $\lceil\sqrt{n}\,\rceil$ colours is established. This result then gives a degree condition for independent coverings in the special case where the graph has $n^2$ vertices. Additionally, this result partially answers Conjecture \ref{Conjecture: Trees} regarding the orthogonal chromatic number of tree graphs.

Orthogonal colourings of tree graphs are interesting because there are only two possible values that the orthogonal chromatic number can take on. Caro and Yuster \cite{caro1999orthogonal} showed that if $T$ is a tree graph with $n$ vertices, then $O\chi(T)=\lceil\sqrt{n}\,\rceil$ or $O\chi(T)=\lceil\sqrt{n}\,\rceil+1$. They proposed Conjecture \ref{Conjecture: Trees} because of an incorrect categorization of the orthogonal number of double star tree graphs. These graphs are now defined and the correct orthogonal chromatic number is determined.

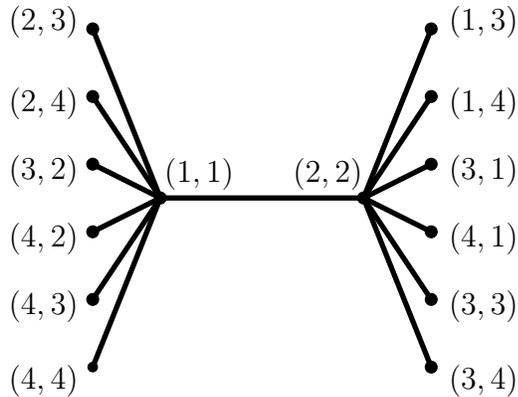
\begin{figure}[h!]
\centering
\begin{tikzpicture}[scale=.9,line cap=round,line join=round,>=triangle 45,x=1cm,y=1cm]
\draw [line width=2pt] (0,5)-- (1,2.5);
\draw [line width=2pt] (1,2.5)-- (0,4);
\draw [line width=2pt] (0,3)-- (1,2.5);
\draw [line width=2pt] (1,2.5)-- (0,2);
\draw [line width=2pt] (1,2.5)-- (0,1);
\draw [line width=2pt] (1,2.5)-- (0,0);
\draw [line width=2pt] (1,2.5)-- (4,2.5);
\draw [line width=2pt] (4,2.5)-- (5,5);
\draw [line width=2pt] (5,4)-- (4,2.5);
\draw [line width=2pt] (4,2.5)-- (5,3);
\draw [line width=2pt] (4,2.5)-- (5,2);
\draw [line width=2pt] (4,2.5)-- (5,1);
\draw [line width=2pt] (4,2.5)-- (5,0);
\draw (-1.4,5.5) node[anchor=north west] {$(2,3)$};
\draw (-1.4,4.3) node[anchor=north west] {$(2,4)$};
\draw (-1.4,3.3) node[anchor=north west] {$(3,2)$};
\draw (-1.4,2.3) node[anchor=north west] {$(4,2)$};
\draw (-1.4,1.3) node[anchor=north west] {$(4,3)$};
\draw (-1.4,.2) node[anchor=north west] {$(4,4)$};
\draw (.9,3.25) node[anchor=north west] {$(1,1)$};
\draw (2.8,3.25) node[anchor=north west] {$(2,2)$};
\draw (5.1,5.5) node[anchor=north west] {$(1,3)$};
\draw (5.1,4.3) node[anchor=north west] {$(1,4)$};
\draw (5.1,3.3) node[anchor=north west] {$(3,1)$};
\draw (5.1,2.3) node[anchor=north west] {$(4,1)$};
\draw (5.1,1.3) node[anchor=north west] {$(3,3)$};
\draw (5.1,.2) node[anchor=north west] {$(3,4)$};
\begin{scriptsize}
\draw [fill=black] (0,0) circle (2pt);
\draw [fill=black] (0,1) circle (2.5pt);
\draw [fill=black] (0,2) circle (2.5pt);
\draw [fill=black] (0,3) circle (2.5pt);
\draw [fill=black] (0,4) circle (2.5pt);
\draw [fill=black] (0,5) circle (2.5pt);
\draw [fill=black] (1,2.5) circle (2.5pt);
\draw [fill=black] (4,2.5) circle (2.5pt);
\draw [fill=black] (5,5) circle (2.5pt);
\draw [fill=black] (5,4) circle (2.5pt);
\draw [fill=black] (5,3) circle (2.5pt);
\draw [fill=black] (5,2) circle (2.5pt);
\draw [fill=black] (5,1) circle (2.5pt);
\draw [fill=black] (5,0) circle (2.5pt);
\end{scriptsize}
\end{tikzpicture}
\caption{Orthogonal Colouring of $D_{14}$}
\label{Figure: Orthogonal Colouring of D14}
\end{figure}

For even $m$, let $D_{m}$ denote the graph obtained by joining the roots of two $K_{1,\frac{m}{2}-1}$ graphs. Caro and Yuster assert in \cite{caro1999orthogonal} that $O\chi(D_{m})=\lceil\sqrt{m}\,\rceil+1$ if $m$ is even satisfying $\lceil\sqrt{m}\,\rceil\lceil\sqrt{m}-1\rceil<m$. The flaw in their proof was that they assumed that no colour could appear $c$ times on leaves, where $c$ is the total number of colours used. This assumption is incorrect, as shown in the following example. For $m=14$, the condition holds, but $O\chi(D_{14})=4$ as shown in Figure \ref{Figure: Orthogonal Colouring of D14}. Also, the colour 4 is used on 4 leaves. The following theorem correctly establishes the orthogonal chromatic number of $D_m$.

\begin{thrm}\label{Theorem: Double Star}
For even $m$, $O\chi(D_m)=\lceil\sqrt{m}\,\rceil=N$ if and only if $m<N^2-1$.
\end{thrm}
\begin{pf}
First, it is shown that for $m\leq N^2-2$, that $O\chi(D_m)=N$. Then, it is shown that for $m\geq N^2-1$, that no orthogonal colouring using $N$ colours exists. In the following, let $x_0$ and $y_0$ denote the root vertices of $D_m$. Then, let $x_1,x_2,\dots,x_{\frac{m}{2}-1}$ and $y_1,y_2,\dots, y_{\frac{m}{2}-1}$ denote the leaves adjacent to $x_0$ and $y_0$ respectively. 

Suppose $m=N^2-2$ in the case where $N$ is even and suppose $m=N^2-3$ in the case where $N$ is odd. In both cases, assign the colour pair $(1,1)$ to $x_0$ and the colour pair $(2,2)$ to $y_0$. This assignment is now extended to the leaves. For $1\leq i \leq N-2$, assign the colour pair $(2,i+2)$ to $x_i$ and the colour pair $(1,i+2)$ to $y_i$. Then for $N-1\leq j \leq 2N-4$, assign the colour pair $(j-N+4,2)$ to $x_j$ and the colour pair $(j-N+4,1)$ to $y_j$. Note that $i+2$ and $j-N+4$ are greater than 2, thus the assigned colour pairs will not cause any colour conflict.

For $3\leq r,s\leq N$, the colour pairs $(r,s)$ can be assigned in any order to the remaining leaves since these pairs will not conflict with the roots $x_0$ and $y_0$. Therefore, in the case where $N$ is even, by arbitrarily assigning all of these colour pairs to the remaining leaves, an orthogonal colouring of $D_m$ using $N$ colours has been constructed. Similarly, in the case where $N$ is odd, by arbitrarily assigning all but one of these colour pairs to the remaining leaves, an orthogonal colouring of $D_m$ using $N$ colours has been constructed.

It is now shown that for $m=N^2$ when $N$ is even, and $m=N^2-1$ when $N$ is odd, that there are no orthogonal colourings using $N$ colours. Let $c_1$ and $c_2$ be two colourings of $D_m$. In $c_1$ (similarly in $c_2$), $x_0$ and $y_0$ must receive different colours. Give $x_0$ the colour pair $(a,b)$ and give $y_0$ the colour $(c,d)$. Then the colour pair $(c,b)$ (similarly $(a,d)$) can not be assigned to any leaf. 

However, in the case where $N$ is even, every colour pair must be used since there are $N^2$ vertices. Thus, no orthogonal colouring using $N$ colours exists. Similarly, in the case where $N$ is odd, all but one colour pair must be used since there are $N^2-1$ vertices, Thus, no orthogonal colouring using $N$ colours exists in this case either. 

Note that if $m< N^2-2$ and $N$ is even, then the orthogonal colouring of $D_{N^2-2}$ constructed in this proof can be restricted to $D_m$ to give an orthogonal colouring using $N$ colours. Similarly, if $m<N^2-3$ and $N$ is odd, then the orthogonal colouring of $D_{N^2-3}$ constructed in this proof can be restricted to $D_m$ to give an orthogonal colouring using $N$ colours. Therefore, for even $m$, $O\chi(D_m)=\lceil\sqrt{m}\,\rceil$ if and only if $m<N^2-1$.
\end{pf}

Theorem \ref{Theorem: Double Star} shows that for some even $m$, there are trees with maximum degree $\frac{m}{2}$ that require $\lceil\sqrt{m}\,\rceil$ colours and also trees that require $\lceil\sqrt{m}\,\rceil+1$ colours. This shows that the maximum degree can not be used to completely classify the orthogonal chromatic number of tree graphs. However, Conjecture \ref{Conjecture: Trees} says that if a tree graph with $n$ vertices has maximum degree less than $\frac{n}{2}$, then an orthogonal colouring using $\lceil\sqrt{n}\,\rceil$ colours exists. This conjecture is false, as shown with the following proposition.

\begin{figure}[h!]
\centering
\begin{tikzpicture}[scale=.7,line cap=round,line join=round,>=triangle 45,x=1cm,y=1cm];
\draw [line width=2pt] (0,0)-- (1,2);
\draw [line width=2pt] (0,2)-- (1,2);
\draw [line width=2pt] (0,4)-- (1,2);
\draw [line width=2pt] (1,2)-- (2,2);
\draw [line width=2pt] (2,2)-- (3,2);
\draw [line width=2pt] (3,2)-- (4,0);
\draw [line width=2pt] (3,2)-- (4,2);
\draw [line width=2pt] (3,2)-- (4,4);
\draw (.9,2) node[anchor=north west] {$u$};
\draw (2.4,2) node[anchor=north west] {$v$};
\begin{scriptsize}
\draw [fill=black] (0,0) circle (2.5pt);
\draw [fill=black] (0,2) circle (2.5pt);
\draw [fill=black] (0,4) circle (2.5pt);
\draw [fill=black] (1,2) circle (2.5pt);
\draw [fill=black] (2,2) circle (2.5pt);
\draw [fill=black] (3,2) circle (2.5pt);
\draw [fill=black] (4,4) circle (2.5pt);
\draw [fill=black] (4,2) circle (2.5pt);
\draw [fill=black] (4,0) circle (2.5pt);
\end{scriptsize}
\end{tikzpicture}
\caption{Counterexample Tree}
\label{Figure: Counter Example}
\end{figure}
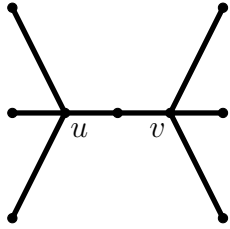

\begin{prop}\label{Proposition: Tree}
For each odd $n$, there exists a tree graph $T$ with $n^2$ vertices such that $\Delta(T)<\frac{n^2}{2}$ and $O\chi(T)=n+1$.
\end{prop}
\begin{pf}
Let $T$ be the tree graph obtained by taking the double star graph $D_{n^2-1}$ and adding a vertex on the edge between the two root vertices. For example, the tree graph obtained for $n=3$ is shown in Figure \ref{Figure: Counter Example}. Suppose that the vertex $u$ has the colour pair $(a,b)$ and the vertex $v$ has the colour pair $(c,d)$ where $(a,b)\neq (c,d)$. Since every other vertex is adjacent to either $u$ or $v$, the colour pair $(a,d)$ and the colour pair $(c,b)$ can not be assigned to a vertex without resulting in a colour conflict. But since the tree graph has $n^2$ vertices, an orthogonal colouring with $n$ colours requires that every colour pair is used. Therefore, an orthogonal colouring of Figure \ref{Figure: Counter Example} with $n$ colours does not exist, and thus $O\chi(T)=n+1$.
\end{pf}

Since the maximum degree of the graph in Figure \ref{Proposition: Tree} is $\lfloor\frac{n}{2}\rfloor<\frac{n}{2}$ for odd $n$, it follows that Conjecture \ref{Conjecture: Trees} is false. This raises the question, if the maximum degree is sufficiently small, does an orthogonal colouring using $\lceil\sqrt{n}\,\rceil$ colours exist? Caro and Yuster \cite{caro1999orthogonal} showed that for any graph $G$, if $\Delta(G)\leq \frac{\sqrt{n}-1}{4}$, then $O\chi(G)=\lceil\sqrt{n}\,\rceil$. The following theorem improves upon this result for tree graphs by using degenerate orderings, which are now defined.

A \textit{$d$-degenerate graph} is a graph such that there exists an ordering of the vertices in which each vertex has $d$ or fewer neighbours that are earlier in the ordering. Such an ordering of the vertices is called a \textit{degenerate ordering}. The \textit{degeneracy} of a graph is the smallest $d$ for which the graph is $d$-degenerate. In particular, tree graphs are known to be $1$-degenerate graphs. The following theorem generalizes a result of Caro and Yuster \cite{caro1999orthogonal}. The same argument is used, except now the degeneracy is used to give a better bound on the orthogonal chromatic number.

\begin{thrm}\label{Theorem: Optimal}
If $G$ is $d$-degenerate with $\Delta(G)<\frac{\sqrt{n}-2d-1}{2}$, then $O\chi(G)=\lceil\sqrt{n}\,\rceil$.
\end{thrm}
\begin{pf}
Consider a $d$-degenerate ordering of the vertices $\{v_1,v_2,\dots,v_n\}$. Let $G_t$ be the graph where all the edges to the vertices $v_{t+1},\dots,v_n$ are removed. Our strategy is to inductively colour $G_t$ with $\lceil\sqrt{n}\,\rceil$ colours, thus colouring $G_n=G$. For $t=1$, $G_1=I_n$, which can be orthogonally coloured with $\lceil\sqrt{n}\,\rceil$ colours \cite{caro1999orthogonal}. 

Now, suppose for $t\geq 1$ that we have a proper orthogonal colouring of $G_{t-1}$ using $\lceil\sqrt{n}\,\rceil$ colours. Assign the same colouring to $G_{t}$. This colouring may not be proper however, because the edges incident to $v_{t}$ from $\{v_1,v_2,\dots, v_{t-1}\}$ are now present. So we show that if this is the case, that we can correct it and still maintain orthogonality.

Let $N_{t}(v_{t})$ be the neighbourhood of $v_{t}$ in $G_{t}$. Then, $|N_{t}(v_{t})|\leq d$, because $v_{t}$ is adjacent to at most $d$ vertices in $\{v_1,v_2,\dots v_{t-1}\}$. Next, let $W$ be the set of vertices $w\in V(G)$, having the property that for some vertex $v\in N_{t}(v_{t})$, $c_1(v)=c_1(w)$ or $c_2(v)=c_2(w)$. Then, $|W|\leq 2|N_{t}(v_{t})|\lceil\sqrt{n}\,\rceil\leq 2d\lceil\sqrt{n}\,\rceil$ because there are $|N_{t}(v_{t})|$ neighbours of $v_{t}$, there are two colourings, and there are at most $\lceil\sqrt{n}\,\rceil$ colours.

Next, let $Y_t$ denote the set of vertices $y\in V(G)\backslash\{v_{t}\}$, having the property that $c_1(y)=c_1(v_{t})$ or $c_2(y)=c_2(v_{t})$. Then, $|Y_t|\leq 2(\lceil\sqrt{n}\,\rceil-1)$ because there are two colourings and there are at most $\lceil\sqrt{n}\,\rceil-1$ colours. Now, let $N(Y_t)$ be the union of open neighbourhoods of these vertices in $G$. Then, $|N(Y_t)|\leq |Y_t|\Delta(G)$. Lastly, let $X=V(G)\backslash(W \cup N(Y_t))$. So $X$ is the set of vertices that do not conflict with the colour assigned to vertices in $N_{t}(v_{t})$ and are not adjacent to vertices that have the same colour as $v_{t}$. The goal is to show $X$ is non-empty. First, note that
\begin{align*}
\Delta(G)& < \frac{\sqrt{n}-2d-1}{2} &&\text{by assumption of the theorem.}\\ 
      & < \frac{\frac{n}{\lceil\sqrt{n}\,\rceil}-2d}{2} &&\text{because $n> \lceil\sqrt{n}\,\rceil\sqrt{n}-\lceil\sqrt{n}\,\rceil$.}\\
      & =\frac{n-2d\lceil\sqrt{n}\,\rceil}{2\lceil\sqrt{n}\,\rceil}   &&\text{by factoring $1/\lceil\sqrt{n}\,\rceil$.}\\
      & <\frac{n-2d\lceil\sqrt{n}\,\rceil}{2\lceil\sqrt{n}\,\rceil-2}  &&\text{because $2\lceil\sqrt{n}\,\rceil > 2\lceil\sqrt{n}\,\rceil -2$.}
\end{align*}
Therefore, we now get the following chain of inequalities:
\begin{align*}
|X|& \geq  n-|W|-|N(Y_t)| \\
   & \geq  n-2d\lceil\sqrt{n}\,\rceil -|Y_t|\Delta(G)\\
   & >  n-2d\lceil\sqrt{n}\,\rceil -(2\lceil\sqrt{n}\,\rceil-2)\left(\frac{n-2d\lceil\sqrt{n}\,\rceil}{2\lceil\sqrt{n}\,\rceil-2}\right)\\
   &=0
\end{align*}

Therefore, the set $X$ is non-empty, so let $x\in X$. Since $x\not\in W$, the colour pair originally assigned to $x$ do not conflict with the neighbours of $v_{t}$. Then, since $x\not\in N(Y_t)$, the colour pair originally assigned to $v_t$ does not conflict with the neighbours of $x$. Therefore, we can interchange the colour pair assigned to $x$ with the colour pair assigned to $v_t$, and still have a proper orthogonal colouring. Thus, we have orthogonally coloured $G_{t}$ with $\lceil\sqrt{n}\,\rceil$ colours. So by induction, $O\chi(G)=\lceil\sqrt{n}\,\rceil$.
\end{pf}

In particular, since tree graphs are $1$-degenerate, by Theorem \ref{Theorem: Optimal}, if $\Delta(T_n)<\frac{\sqrt{n}-3}{2}$, then $O\chi(T_n)=\lceil\sqrt{n}\,\rceil$. Note that in the case where the graph has $n^2$ vertices, the orthogonal colouring created corresponds to an independent covering since all colour classes have the same size. It remains an open problem to determine if $\frac{\sqrt{n}-3}{2}$ is the best possible upper bound for tree graphs. 

\bibliographystyle{amsplain}

\bibliography{Independent_Coverings_and_Orthogonal_Colourings_References}

\end{document}